\documentclass[11pt]{article}

\usepackage{amsmath}
\usepackage{amssymb}
\usepackage{amscd}

\newtheorem{theo}{Theorem}[section]

\newtheorem{remar}[theo]{Remark}
\newtheorem{prop}[theo]{Proposition}

\newtheorem{lemma}[theo]{Lemma}

\newtheorem{Example}[theo]{Example}

\newcommand {\wt}{\widetilde}

\newcommand{\Rea}{\operatorname{Re}}

\newcommand{\fdim}{\hspace*{\fill}$\Box$}
\newcommand{\dimostr}{{\bf Proof: }}

\newcommand{\R}{\Bbb{R}}

\newcommand{\complex}{\Bbb{C}}
\newcommand{\C}{\Bbb{C}}

\newcommand{\N}{\Bbb{N}}

\newcommand{\Tanh}{{\rm Tanh}}
\newcommand{\Cosh}{{\rm Cosh}}
\newcommand{\K}{K\"{a}hler}

\begin{document}

\title{Riemannian Geometry of Hartogs domains
\footnote{This work was supported by the MIUR Project \lq\lq
Riemannian metrics and differentiable manifolds'' .}}

\author{Antonio J. Di Scala,  Andrea Loi, Fabio Zuddas} \date{}
\maketitle

\abstract
Let
$D_F = \{ (z_0, z) \in {\C}^{n} \; | \; |z_0|^2 < b, \|z\|^2  < F(|z_0|^2)
\}$ be  a strongly pseudoconvex  Hartogs domain
endowed with the \K\ metric $g_F$ associated to the \K\ form
$\omega_F = -\frac{i}{2} \partial \overline{\partial} \log
\left(F(|z_0|^2) - \|z\|^2\right)$.

This paper  contains several results on the Riemannian geometry
of these domains.
These are  summarized in   Theorem \ref{mainteor},
Theorem  \ref{mainteor1} and Theorem  \ref{mainteor2}.
 In the
first one we prove  that if $D_F$ admits a non special   geodesic
(see definition below) through the origin whose trace is a
straight line then $D_F$ is holomorphically isometric to an open subset of  the
complex hyperbolic space. In the second theorem  we prove that
all the  geodesics through the origin of $D_F$  do not self-intersect,
we find necessary and sufficient
conditions on $F$ for $D_F$ to be geodesically complete and
we prove that $D_F$ is locally irreducible as a Riemannian manifold.
Finally, in Theorem \ref{mainteor2},
 we compare the Bergman metric $g_B$ and the metric $g_F$
in a bounded Hartogs domain and we prove that if $g_B$
is a multiple of $g_F$, namely
$g_B=\lambda g_F$,  for some  $\lambda\in \R^+$,  then $D_F$ is holomorphically isometric to
an open subset of  the
complex hyperbolic space.

\vskip 0.1cm


\noindent
{\it{Subj.Class}}: 53C55, 32Q15, 32T15.

\noindent

\section{Introduction}
Let $b \in \R^+ \cup \{ + \infty \}$ and let $F: [0,
b) \rightarrow (0, + \infty)$ be a non increasing smooth
function  on $[0, b)$. The $n$-dimensional  {\em Hartogs domain}
$D_F\subset {\C}^{n}$ associated to the function $F$ is defined by
\begin{equation}\label{hartogs}
D_F = \{ (z_0, z) \in {\C}^{n} \; | \; |z_0|^2 < b, \|z\|^2  < F(|z_0|^2)
\},
\end{equation}
where $z=(z_1,\dots ,z_{n-1})$ and $\|z\|^2= |z_1|^2 +\cdots+ |z_{n-1}|^2$.

Under the assumption that $D_F$ is strongly pseudoconvex one can
prove (see Proposition 2.1 in \cite{secondoartic}) that the
natural $(1, 1)$-form on $D_F$ given by
\begin{equation}\label{omegaf}
\omega_F = -\frac{i}{2} \partial \overline{\partial} \log
\left(F(|z_0|^2) - \|z\|^2\right)
\end{equation}
is  a \K\ form on $D_F$ and this is equivalent to the requirement
that  $F$ satisfies the condition
\begin{equation}\label{Kcond}
\left(\frac{xF'}{F}\right )' <0, \ \ \   x \in [0, b),
\end{equation}
where the derivatives are taken  with respect to the variable
$x=|z_0|^2$.
Notice that (\ref{Kcond}) (and hence the strongly pseudoconvexity of $D_F$)
turns out to be equivalent to the
strongly pseudoconvexity of  the boundary of $D_F$
at all $z = (z_0, z_1,\dots,z_{n-1})$ with
$|z_0|^2 < x_0$ (see Proposition 2.1 in \cite{secondoartic} for a proof).
Denote by  $g_F$ the \K\ metric associated to the \K\
form $\omega_F$. Throughout all this paper we assume that $D_F$ is
equipped with this \K\ metric. Notice that when $F(x)=1-x, 0\leq
x< 1$, $D_F$ equals the $n$-dimensional complex hyperbolic space
${\complex} H^n$, namely the unit ball $B^n$ in ${\C}^n$ equipped
with the hyperbolic metric $g_{hyp}=g_F$. Hartogs domains are
interesting both from the mathematical and the physical point of
view (see for example \cite{me}, \cite{constscal} and
\cite{secondoartic} for the study of some Riemannian properties of
$g_F$ and the Berezin quantization of $(D_F, g_F)$, \cite{alcu}
and \cite{primoartic} for the construction of global symplectic
coordinates on these domains and \cite{roos} for the construction
of \K\--Einstein metrics on Hartogs type domains on symmetric
spaces).

\vskip 0,3cm

In this paper we are interested in  the Riemannian geometry of Hartogs domains.
In particular we study the geodesics,  the completeness and local irreducibility
 of such domains with respect to the metric $g_F$.
We denote by  ${\cal G}$ the set of  those geodesics passing through the origin
whose traces are   straight lines of ${\complex}^n$ intersected with $D_F$.
Since the isometry group of $D_F$
contains $U(1)\times U(n-1)$,
it  is easily seen that
the set ${\cal S}$ of  geodesics of $D_F$
passing through the origin and
contained   in the plane  $z_0=0$
or in the complex line
$z_1=\cdots =z_{n-1}=0$
is included in   ${\cal G}$.
 A geodesic $\ell \in {\cal G}$ will be called a
{\em special  geodesic} if it  belongs  to  ${\cal S}$.

Our first  result is the following interesting characterization
of the complex  hyperbolic space amongst Hartogs domains.
\begin{theo}\label{mainteor}
Let $(D_F, g_F)$ be a  Hartogs domain.
If  there exists $\ell\in {\cal G}$ such that $\ell\notin {\cal S}$,
then $D_F$ is holomorphically isometric to an open subset of
${\C} H^n$.
\end{theo}
In other words the previous theorem  asserts that  if there exists
one non special  geodesic  $\ell$ through the origin of
$D_F$ whose trace is  a straight line, then $D_F\subset{\complex}H^n$
(and hence, in this case,  ${\cal G}$ coincide with the set of all the geodesic through the origin).

Our second result is the following:
\begin{theo}\label{mainteor1}
Let $(D_F, g_F)$ be a Hartogs domain.
Then the following facts hold true.
\begin{itemize}
\item [(i)]
All the  geodesics through the origin of $D_F$  do not self-intersect;
\item [(ii)]
$D_F$ is geodesically complete
with respect to the \K \ metric $g_F$
if and only if
\begin{equation}\label{complcond}
\int_0^{\sqrt{b}} \sqrt{- \left(\frac{xF'}{F}\right )'}|_{x=
u^2} \ du = +\infty,
\end{equation}
where we define  $\sqrt{b}=+\infty$ for $b=+\infty$;
\item [(iii)]
$(D_F, g_F)$ is locally irreducible around any of its  points.
\end{itemize}
\end{theo}

The first part of the previous theorem  should be compared with the main result of D'Atri and Zhao
\cite{datri} asserting that  in
 a bounded homogeneous
domain equipped with its Bergman metric
all  the  geodesics  do not intersect.
(Notice that the homogeneous assumption for an Hartogs domain implies
that $D_F$ is holomorphically equivalent to $B^n$,
 the unit ball in ${\complex}^n$ (see e.g. Theorem 6.11 in \cite{ik}).

Other properties of the geodesics of the Bergman metric can be
found in \cite{feff} and \cite{herb}. In \cite{feff} Fefferman
deeply studied the geodesics of the Bergman metric at the boundary
points of  a bounded domain $D$ while \cite{herb} is concerned
with   the existence of a closed geodesic in any non trivial
homotopy class of a (non simply-connected) bounded domain.
Regarding the completeness of the Bergman metric on a bounded domain,
 the reader is referred to the classical paper of S. Kobayashi \cite{ko}.

\vskip 0.3cm

By the previous discussion  it is natural to compare the Bergman metric $g_B$
and the metric $g_F$ on a {\em bounded} Hartogs domain.
Our third and last result is the following:

\begin{theo}\label{mainteor2}
Let $D_F$ be a bounded Hartogs domain.
Assume that $g_B$ is a multiple of $g_F$, namely
$g_B=\lambda g_F$,
for some $\lambda\in\R^+$. Then
$g_F$ is \K --Einstein and therefore
$D_F$ is holomorphically isometric to
an open subset of ${\C}H^n$.
\end{theo}

The first part of the proof of the previous theorem is an adaptation of the proof of the  following   (unpublished) proposition  communicated by  Miroslav Engli\v{s}  to the second author and which deals
with the more general  class of  {\em generalized} Hartogs domains.

\begin{prop}(Engli\v{s}) \label{EnglisG} Let
$$ \wt \Omega = \{(z,w)\in\Omega\times\mathbb{C}^k: \|w\|^2 <F(z)\}  $$
be a bounded and simply-connected
 generalized Hartogs domain,
where $\Omega$  is
a  pseudoconvex
domain in $\mathbb{C}^n$ and $-\log F$
is a smooth strictly-PSH function
on $\Omega$.
Let $g_B$ be the Bergman metric and
let $g_F$ be  the \K\  metric on
 $\wt\Omega$ whose \K\ potential is
  $-\log (F(z) - \|w\|^2)$.
If  $g_B=\lambda g_F$, for some $\lambda\in {\R}^+$,
then $g_F$   is K\"ahler-Einstein.
\end{prop}

The next section is dedicated to the proof of our theorems.

\section{Proof of the main results}
The following lemma is the  main tool in the proofs of Theorem \ref{mainteor}
and Theorem \ref{mainteor1}.
\begin{lemma}\label{mainlemma}
Let $(D_F, g_F)$
be a Hartogs domain.
 Let
$M\subset D_F$ be the real (plane)
surface given by:
\begin{equation}\label{m}
M=D_F\cap
\{Im (z_0)=Im (z_1)=0, \, z_j=0, j=2, \dots , n-1\},
\end{equation}
and denote by $g$ the metric  induced  on $M$ by $g_F$.
Then $(M, g)$ is totally geodesic, has constant Gaussian curvature
equal to
$-\frac{1}{2}$ and is geodesically complete if and only if
condition (\ref{complcond}) above  is satisfied.
\end{lemma}
\dimostr
The surface $M$ is the fixed point set of the isometry of
$D_F$ given by $(z_0, z_1, z_2, \dots , z_{n-1})\mapsto (\bar{z}_0, \bar{z}_1, -z_2, \dots , -z_{n-1})$
and hence it is totally geodesic in $D_F$.
By  setting $u=Re (z_0)$ and $v=Re  (z_1)$, this surface can be described as
\begin{equation}\label{superficieM}
M=\{(u, v )\in {\R}^2|\ v^2<F(u^2), u^2<b\}.
\end{equation}
Furthermore, it is not difficult to see
that the metric $g$ induced by  $g_F$
on $M$ is given by
\begin{equation}\label{matrixreal2}
g= \left( \begin{array}{cc}
g_{11} & g_{12}\\
g_{21}  & g_{22}
\end{array} \right)
= \frac{2}{(F - v^2)^2} \left( \begin{array}{cc}
C & - F' u v\\
- F' u v  & F
\end{array} \right),
\end{equation}
where $C = F'^2 \cdot u^2 - (F' + F'' \cdot u^2)(F - v^2 )$
and $F$
and its derivatives are evaluated at $u^2$.
By a straightforward,  but  long computation,
one can verify that the Gaussian curvature of $g$ equals    $-1/2$.
Hence $(M, g)$ is
isometric to an open subset, say $U$, of
${\R}H^2(\tiny{-\frac{1}{2}})$,
namely  the unit  disk
$\{(x,y) \ | \ x^2 + y^2 < 1 \}$ in ${\R}^2$
endowed with the Beltrami-Klein metric
\begin{equation}\label{bk}
g_{BK}=\frac{2}{(1 - x^2 - y^2)^2}\left[ (1-y^2) dx^2 + 2 x ydx dy +(1-x^2)dy^2\right].
\end{equation}
An  isometry between $(M, g)$ and $U$ can be described explicitly.
Indeed, let
$\psi: (- \sqrt{b}, \sqrt{b}) \rightarrow
{\R}$ be the strictly  increasing real valued  function defined by
$$\psi(u) = \int_0^u \sqrt{-  \left( \frac{x F'}{F}\
\right)'}|_{x=s^2} ds .$$
Then, it is not hard to see  that the map
$$\Psi: M \rightarrow {\R}H^2(\tiny{- \frac{1}{2}}),
(u,v) \mapsto \left( \Tanh(\psi(u)), \frac{v}{\Cosh(\psi(u))
\sqrt{F(u^2)}} \right)$$
is an injective local diffeomorphism satisfying
$\Psi^*(g_{BK})=g$.
Therefore,  the completeness of  $M$ is equivalent to  $\Psi(M)={\R}H^2$,
which is easily seen to be  equivalent
to  condition
(\ref{complcond}),
and we are done.
\fdim

\begin{remar}\label{remarmainteor}\rm
The fact  that the surface $M$ in Lemma  \ref{mainlemma}
is totally geodesic and that the isometry group of $D_F$
contains $U(1)\times U(n-1)$ implies the existence of an isometry
of $D_F$, fixing the origin and taking any  given geodesic passing through the origin
of $D_F$ to a geodesic lying in $M$.
This will be a key point  in the proofs of both Theorem \ref{mainteor}
and Theorem \ref{mainteor1}.
\end{remar}

\begin{remar} \label{Mzero} \rm
All the $n$-dimensional Hartogs
domain $D_F$ contains  the complex totally
geodesic surface
\begin{equation}\label{m0}
{\cal B} = \{ z_j = 0, \ j=1,\dots ,n-1  \} \cap D_F,
\end{equation}
which in the literature of complex analysis  is called the {\em base}
 of the Hartogs domain $D_F$.
In view of the previous lemma,
it  is natural to consider  the  Hartogs domains where the
Gaussian
curvature of ${\cal B}$ is constant, ,  say equal to $K_0$,.
It is not hard to see that such   domains can be classified
as follows:
\begin{itemize}
\item [(a)] if $K_0=0$ then  $F(t)=ce^{-kt}, c, k>0$, $t\in [0,
+\infty )$, (complex analysts often refer to these
domains as  {\em Spring domains});
 \item [(b)] If $K_0>0$ then $F(t) = (c_1 + c_2
t)^{- \frac{2}{ K_0}}$, with $c_1
>0$, $c_2 >0$, $t\in [0, +\infty )$;
\item [(c)] If   $K_0 <0$ then $F(t) = (c_1 + c_2
t)^{- \frac{2}{ K_0}}$,
with $c_1 >0$, $c_2 <0$, $t\in [0, -\frac{c_1}{c_2} )$.
\end{itemize}

Notice that in the case (b),
the corresponding Hartogs domain $D_F$
 cannot  be geodesically complete.
 In fact in this case also its base
${\cal B}$
 would be complete and hence biholomorphic to ${\C}P^1$,
 yielding the contradiction
 ${\cal B}\cong{\C}P^1\subset D_F\subset {\C}^n$
 (cfr.  Example \ref{fin} at the end of the paper).
\end{remar}

\vskip 0.3cm

\noindent
{\bf Proof of Theorem \ref{mainteor}. \rm}

\noindent
Let  $\ell$ be  a geodesic
as in the statement of the theorem.
Since $\ell\notin {\cal S}$,
by Remark \ref{remarmainteor},
we can assume $\ell\subset M$
and that
$$\ell =\{v=ku, k\neq 0\}\cap M, $$
where $u$ and $v$ are the parameters
introduced in the proof of Lemma \ref{mainlemma}.
Hence $\ell$ can be parametrized  as
$t\mapsto (u(t), v(t)=ku(t))$, where
$t$ varies in a real interval containing the origin
and  the following  geodesic equations have to be satisfied
$$u'' + \Gamma^1_{11} u'^2 + 2 \Gamma^1_{12} u' v' + \Gamma^1_{22}
v'^2 = 0, \  \
v'' + \Gamma^2_{11} u'^2 + 2 \Gamma^2_{12} u' v' + \Gamma^2_{22}
v'^2 = 0,  $$
namely

\begin{equation}\label{geod1ret}
u'' + \Gamma^1_{11} u'^2 + 2 k \Gamma^1_{12} u'^2 + k^2
\Gamma^1_{22} u'^2 = 0
\end{equation}

\begin{equation}\label{geod2ret}
k u'' + \Gamma^2_{11} u'^2 + 2 k \Gamma^2_{12} u'^2 + k^2
\Gamma^2_{22} u'^2 = 0,
\end{equation}
where $\Gamma^i_{jk}, i, j, k=1, 2$ are the Christoffel symbols
(see e.g. \cite{docarmo}).
A straightforward computation gives :

\begin{footnotesize}

$$\Gamma^1_{11}  =  \frac{1}{2D} \left( g_{22} \frac{\partial
g_{11}}{\partial u} - g_{12}\left(2 \frac{\partial
g_{12}}{\partial u} -\frac{\partial g_{11}}{\partial v}  \right)
\right) =$$

\begin{equation}\label{unounouno}
= \frac{-4 u}{D (v^2 - F)^4}
[ u^2(2
F'^2 + v^2 F'') -F(v^2 - F)(2 F'' + u^2 F''') - F F'(2 F' + 3u^2
F'') ],
\end{equation}

\begin{eqnarray}
\Gamma^2_{11} &  = & \frac{1}{2D} \left( - g_{12} \frac{\partial
g_{11}}{\partial u} + g_{11} \left(2 \frac{\partial
g_{12}}{\partial u} -  \frac{\partial g_{11}}{\partial v} \right)
\right) = \nonumber \\ & = & \frac{4 u^2 v}{D (v^2 - F)^3} [-u^2
F''^2 + F'(F'' + u^2 F''')],
\end{eqnarray}

\begin{eqnarray}
\Gamma^1_{12} & = & \frac{1}{2D} \left(  g_{22} \frac{\partial
g_{11}}{\partial v} -  g_{12} \frac{\partial  g_{22}}{\partial u}
\right) = \nonumber \\ & = & \frac{-4 v}{D (v^2 - F)^4} [-u^2 F'^2
+ F(F' + u^2 F'')],
\end{eqnarray}

\begin{eqnarray}
\Gamma^2_{12} & = & \frac{1}{2D} \left(  g_{11} \frac{\partial
 g_{22}}{\partial u} -  g_{12} \frac{\partial g_{11}}{\partial v}
\right) = \nonumber \\ & = & \frac{4 u F'}{D (v^2 - F)^4} [-u^2
F'^2 + F(F' + u^2 F'')],
\end{eqnarray}

\begin{eqnarray}
\Gamma^1_{22} = \frac{1}{2D} \left( - g_{12}\frac{\partial
g_{22}}{\partial v} +  g_{22} \left( 2\frac{\partial
g_{12}}{\partial v} -  \frac{\partial  g_{22}}{\partial u} \right)
\right) = 0,
\end{eqnarray}

\begin{eqnarray}\label{dueduedue}
\Gamma^2_{22} & = & \frac{1}{2D} \left(  g_{11} \frac{\partial
g_{22}}{\partial v} - g_{12}\left(2 \frac{\partial
g_{12}}{\partial v} -  \frac{\partial  g_{22}}{\partial u} \right)
\right) = \nonumber \\ & = & \frac{-8 v}{D (v^2 - F)^4} [-u^2 F'^2
+ F(F' + u^2 F'')] ,
\end{eqnarray}
\end{footnotesize}
where
$$D = g_{11}g_{22} - g_{12}^2=
4\frac{CF-F'^2u^2v^2}{(F-v^2)^4}.$$
By solving (\ref{geod1ret}) with respect to $u''$ and substituting
into (\ref{geod2ret}) we get
\begin{equation}
u'^2 [ \Gamma^2_{11} + k (2 \Gamma^2_{12} - \Gamma^1_{11} ) + k^2
( \Gamma^2_{22} - 2 \Gamma^1_{12} ) - k^3 \Gamma^1_{22}] = 0
\end{equation}

Since  $u' \neq 0$ we get

\begin{equation}\label{esercizioDoCarmo}
\Gamma^2_{11} + k (2 \Gamma^2_{12} - \Gamma^1_{11} ) + k^2 (
\Gamma^2_{22} - 2 \Gamma^1_{12} ) - k^3 \Gamma^1_{22} = 0
\end{equation}

(where $\Gamma^k_{ij} = \Gamma^k_{ij}(u, ku)$).

By using  (\ref{unounouno}) -  (\ref{dueduedue}),
after a very long but straightforward calculation,
the previous equation becomes
\begin{equation}\label{equazDiffer}
\frac{8 k u\left ( u^4 F''^2 + F (2 F'' + u^2 F''') - F'(2 u^2 F''
+ u^4 F''') \right) }{D (k^2 u^2 - F)^3} = 0,
\end{equation}
which, by setting $u^2 = t, 0\leq t<b$,   is equivalent
to the following  ODE

\begin{equation}\label{equazDiffer2}
t^2 F''^2 + F (2 F'' + t F''') - F'(2 t F'' + t^2 F''') =0.
\end{equation}

Notice  that for $t \neq0$
this equation can be written as

\begin{equation}\label{equazDiffer3}
t^2 F''^2 + \left(\frac{F}{t} - F' \right) (t^2 F'')' =0.
\end{equation}

By setting $G = t^2 F''$
equation (\ref{equazDiffer3}) becomes
\begin{equation}\label{equazDiffer4}
G' = - \frac{F'' t}{F - F' t}  G
\end{equation}
(notice that $F- F' t >0$ for all $0<t <b$
since $F$ is not decreasing)
and hence
\begin{equation}\label{equazDiffer4}
G(t) = c \ e^{\int \frac{-F'' t}{F - F't} dt} = c \ (F - F't),
\end{equation}
 for some
$c \in \R$.
For  $t\rightarrow 0$ we
get $c F(0) = 0$, i.e. $c = 0$.
Therefore   $G = t^2 F'' = 0$,  which
implies
$F (t)= c_1 - c_2 t $ for some $c_1,
c_2 >0$.
Then the map
$$\phi: D_F \rightarrow {\C}H^{n}, \ (z_0, z_1,\dots,z_{n-1}) \mapsto
\left( \frac{z_0}{\sqrt{c_1/c_2}}, \frac{z_1}{\sqrt{c_1}},\dots,\frac{z_{n-1}}{\sqrt{c_1}} \right) $$
is a holomorphic isometry  of $D_F$
into an open subset of  ${\C}H^n$ and this  concludes the proof  of Theorem \ref{mainteor}.
\fdim

\begin{remar}\rm
In the very definition of a Hartogs domain $D_F$
we have assumed (cfr. the introduction)
that $F$ is non increasing in the interval
$[0, b)$.
The statement of Theorem \ref{mainteor}
holds true also without this assumption.
Indeed, it follows  by condition (\ref{Kcond})
that $F'(t)<0$ in a suitable interval
 $0\leq t<\epsilon< b$, for some $\epsilon$,
and the proof  works also in this case.
\end{remar}

\vskip 0.3cm

\noindent
{\bf Proof of Theorem \ref{mainteor1}}

\noindent
Let $\ell\subset D_F$ be a geodesic passing through the origin.
By Remark \ref{remarmainteor} we can assume $\ell\subset M$.
On the other hand by Lemma \ref{mainlemma},  $(M, g)$
is isometric to an open subset of ${\R}H^2(-\frac{1}{2})$
where it is well-known  that all the geodesics do not self intersect.
This proves (i) of Theorem \ref{mainteor1}.

Notice that,  again by Remark \ref{remarmainteor}
and by Hopf--Rinow's theorem the completeness of $g_F$
is equivalent to that of $g$, which by Lemma \ref{mainlemma}
is equivalent to  (\ref{complcond}) and we this proves (ii).

In order to prove (iii), assume by contradiction that $D_F$ is
locally reducible around some point, say $p\in D_F$. Since the
group $U(1) \times U(n-1)$ acts by isometries on $(D_F,g_F)$  we
can assume that $p \in M$ where $M$ is given by (\ref{m}). So $p =
(u, v  , 0, \cdots, 0)$ and we can assume that both  $u, v$ are
(real numbers) different from zero (indeed if one of them is zero,
say $p_1=0$, then $D_F$ is locally reducible around the point
$p'=(p_1', p_2, 0, \dots 0)$ with $p_1'$ sufficiently close to
zero). Therefore there exists an neighborhood $D\subset D_F$ of $p
\in D_F$ such that $(D, g_F)$ splits as a Riemannian product i.e.
$D = A \times B$, where $A$ and $B$ are K\"ahler manifolds. So the
Lie algebra $\frak{g}$ of Killing vector fields of $D$ also splits
into two (or more) factors. Since $\frak{u}(1) \times
\frak{u}(n-1) \subset \frak{g}$ it follows that $\frak{g}$ has at
most two factors. Moreover since $p = (u, v  , 0, \cdots, 0)$ with
$u,v \neq 0$ we can recover the tangent space to the Riemannian
factors $A$ and $B$. Thus, the factor $A$ is an open subset $A
\subset {\C}$, with $u\in A$, and $B$ is an open subset $B\subset
{\C}^{n-1}$, with $(v , 0, \cdots , 0)\in B$. In particular such
Riemannian factors must be orthogonal w.r. to $g_F$. Then the
coefficient $g_{12}$ of the metric $g$ on $M$ induced by $g_F$ has
to be zero for $u$ and $v$ different from zero. On the other hand,
by (\ref{matrixreal2}) above, $g_{12}=- F'uv\neq 0$, a
contradiction. This concludes the proof of (iii). \fdim

\vskip 0.2cm

\noindent
With  the aid of  (ii) in Theorem \ref{mainteor1} we now  study the completeness of
two specific Hartogs domains.

\begin{Example}\rm
If $F(t) = ce^{-kt}$, $c, k>0$, $t \in [0, + \infty)$ then condition (\ref{complcond}) is easily seen to
be satisfied, so we get that the the Spring domains are complete  (cfr.  (a) of Remark \ref{Mzero}).
\end{Example}
\begin{Example}\label{fin}\rm
If $F(t) = \frac{1}{(c_1 + c_2 t)^p}$ ($p \in {\N}^+$),
$t \in [0, + \infty)$ , then
$$\int_0^{\sqrt{b}} \sqrt{- \left(\frac{xF'}{F}\right )'}|_{x=
u^2} \ du = \int_0^{+ \infty} \frac{\sqrt{c_1 c_2 p}}{c_1 + c_2
u^2} du = \frac{\pi}{2} \sqrt{p} < \infty$$
which proves that, for such $F$, the domain $D_F$ is not
complete (cfr.  (b) of Remark \ref{Mzero}).
\end{Example}

\vskip 0.3cm

We now prove the last result of this paper.

\noindent
{\bf Proof of Theorem \ref{mainteor2}}

\noindent
Recall that the Bergman metric
$g_B$ on $D_F$ is, by definition, the one given by the
K\"ahler potential $\log\wt K(z_0,z;z_0,z)$, where $\wt K(z_0,z ; z_0',z')$
is the Bergman kernel of $D_F$.
Let
\begin{equation}\label{tildeF}
\tilde F(z_0, z) := F(|z_0|^2) -\|z\|^2 .
\end{equation}
Note that this is a local defining function (positively signed)
for $D_F$ at any boundary point $(z_0,z)$ with $z\in {\cal B}$,
and such boundary points are strictly pseudoconvex.
The hypothesis of the theorem  and the fact that $D_F$ is contractible means that
$$ \log \wt K(z_0,z) = -\lambda \log \wt F(z_0,z) + 2 \Rea G(z_0,z)  $$
for some  holomorphic function $G$ on $D_F$; here and below we will
write just $\wt K(z_0,z)$ for $\wt K(z_0,z; z_0,z)$. By rotational
symmetry of $\wt K$ and $F$, the pluriharmonic function $2\Rea G$
must depend only on $|z_0|^2$ and $\|z\|^2$, hence must be a positive  constant, say $\mu$.
Thus
\begin{equation}\label{tildeK}
 \wt K(z_0,z) = \frac{\mu }{\wt F(z_0,z)^{\lambda}}.
 \end{equation}
On the other hand, by Fefferman's formula \cite{feff} for the boundary
singularity of the Bergman kernel,
\begin{equation}\label{fefferman}
 \wt K(z_0,z) = \frac{a(z_0,z)}{\wt F(z_0,z)^{n+1}} + b(z_0,z) \log \wt F(z_0,z),
 \qquad (z_0,z)\in D_F,
 \end{equation}
where $a,b\in C^\infty({\cal B}\times \mathbb{C}^{n-1})$ and
\begin{equation}\label{az0}
a(z_0,z) = \frac{n!}{\pi^{n}} \; J[\tilde F](z_0, z),
 \end{equation}
 for  $z_0\in {\cal B}$ and  $\|z\|^2=F(|z_0|^2)$
 and
where $J[\wt F]$ is the Monge-Ampere determinant
\[ J[\wt F] = (-1)^{n} \det\bmatrix\wt  F &\frac{\partial \wt  F}{ \partial z_0} & \partial_{z} \wt F \\
\frac{\partial \wt  F} {\partial\bar z_0}&\frac{\partial^2 \wt  F} {\partial z_0  \partial\bar z_0} &
 \partial_{z}(\frac{\partial\wt F}{\partial\bar z_0}) \\
\partial_{\bar z} \wt  F &  \partial_{\bar z}(\frac{\partial\wt F}{\partial\bar z_0}) & \partial_{\bar z} \partial_{z} \wt F
\endbmatrix. \]
A direct  computation gives
 \begin{equation}\label{J}
 J[\wt F] =  -F^{2} \frac{\partial^2\log F}{\partial z_0\partial\bar z_0}.
 \end{equation}
 (which
 depends only on $|z_0|^2$).
By comparing (\ref{tildeK})
with (\ref{fefferman}) one gets:
$$\mu = \frac{a(z_0, z) \wt F(z_0,z)^{\lambda}}{\wt F(z_0,z)^{n+1}} + b(z_0,z)\wt F(z_0,z)^{\lambda} \log \wt F(z_0,z),  \qquad (z_0,z)\in D_F, $$
which evaluated at
$\|z\|^2=F(|z_0|^2)$,
forces $\lambda =n+1$. Further,  by  (\ref{az0}) and (\ref{J}), the last expression gives
$$-F^2\frac{\partial^2\log F}{\partial z_0\bar\partial z_0}=c, $$ for
all $z_0\in {\cal B}$ and $\|z\|^2=F(|z_0|^2)$,
where $c$ is the negative constant given by $c=-\frac{\mu\pi^n}{n!}$
(notice that  the   condition $\|z\|^2=F(|z_0|^2)$
 is superfluous, since nothing there  depends on $z$).
Feeding this back into  formula (\ref{J}) one gets
$J[F](z_0,z)=c$ for all $(z_0,z)\in D_F$,
i.e. $g_F$ is K\"ahler-Einstein.

\vskip 0.1cm

Let us recall now Lemma 3.1. of \cite{secondoartic}.

\begin{lemma}\label{lemmafin}
 Let $(D_F,g_F)$ be an $n$-dimensional Hartogs
domain. Assume that one of its generalized scalar curvatures is
constant. Then $(D_F,g_F)$ is holomorphically isometric to
an open subset of the
$n$-dimensional complex hyperbolic space.
\end{lemma}

Since the scalar curvature is one of the generalized scalar
curvatures the proof of Theorem \ref{mainteor2} is complete.
\fdim

\vskip 0.3cm

\small{}

\vskip.3cm

\noindent 
Dipartimento di Matematica,
Politecnico di Torino,\\
Corso Duca degli Abruzzi 24, 10129 Torino, Milano, Italy,  \\
\textit{E-mail:} \texttt{antonio.discala@polito.it}

\vskip.3cm

\noindent
Dipartimento di Matematica e Informatica
Universit\`a di Cagliari,\\
Via Ospedale 72, 09124 Cagliari, Italy, \\
\textit{E-mail:} \texttt{loi@unica.it}

\vskip.3cm

\noindent
Dipartimento di Matematica e Informatica
Universit\`a di Cagliari,\\
Via Ospedale 72, 09124 Cagliari, Italy, \\
\textit{E-mail:} \texttt{fzuddas@unica.it}

\end{document}